\documentclass[a4paper]{article}

\usepackage[margin=1in]{geometry} 

\usepackage[T1]{fontenc}
\newcommand{\changefont}[3]{
\fontfamily{#1} \fontseries{#2} \fontshape{#3} \selectfont}

\changefont{ptm}{m}{n}

\usepackage{setspace} 
\usepackage{graphicx}

\usepackage{amsfonts}
\usepackage{amsmath}
\usepackage{amssymb}
\usepackage{graphics}
\usepackage{mathrsfs}
\usepackage{color}
\newtheorem{remark}{Remark}[section]

\newtheorem{theorem}{Theorem}[section]

\newtheorem{lemma}{Lemma}[section]

\newtheorem{definition}{Definition}[section]

\long\def\symbolfootnote[#1]#2{\begingroup%
\def\thefootnote{\fnsymbol{footnote}}\footnote[#1]{#2}\endgroup} 

\begin{document}

\begin{center}
\Large \textbf{Modulo periodic Poisson stable solutions of dynamic equations on a time scale}
\end{center}

\begin{center}
\normalsize \textbf{Fatma Tokmak Fen$^{1}$, Mehmet Onur Fen$^2,$\symbolfootnote[1]{Corresponding Author. E-mail: monur.fen@gmail.com}} \\
\vspace{0.2cm}
\textit{\textbf{$^1$Department of Mathematics, Gazi University, 06560 Ankara, Turkey}} \\

\vspace{0.1cm}
\textit{\textbf{$^2$Department of Mathematics, TED University, 06420 Ankara, Turkey}} \\
\vspace{0.1cm}
\end{center}

\vspace{0.3cm}

\begin{center}
\textbf{Abstract} 
\end{center}

The existence, uniqueness, and asymptotic stability of modulo periodic Poisson stable solutions of dynamic equations on a periodic time scale are investigated. The model under investigation involves a term which is constructed via a Poisson stable sequence. Novel definitions for Poisson stable as well as modulo periodic Poisson stable functions on time scales are provided, and the reduction technique to systems of impulsive differential equations is utilized to achieve the main result. An example which confirms the theoretical results is provided.

\vspace{-0.2cm}

\noindent\ignorespaces

\vspace{0.3cm}
 
\noindent\ignorespaces \textbf{Keywords:} Modulo periodic Poisson stability, Dynamic equations, Impulsive differential equations, Periodic time scale

\vspace{0.2cm}

\noindent\ignorespaces \textbf{Mathematics Subject Classification:} 34N05, 34A37 

\vspace{0.6cm}


\section{Introduction and preliminaries} \label{sec1}

Poisson stable motions, which were first introduced by Poincar\'{e} \cite{Poincare1892}, include the cases of oscillations such as periodic, quasi-periodic, almost periodic, almost automorphic, recurrent, and pseudo-recurrent ones \cite{ Birkhoff66, Bochner64, Corduneanu09, Knight81, Veech63}. Results on Poisson stable solutions for stochastic differential equations and a class of fourth-order dynamical systems can be found in the studies \cite{Cheban20, Liu22, Pchelintsev21}. Recently, a new type of flow called modulo periodic Poisson stable (MPPS) was introduced in paper \cite{Akhmet21}, where the authors also considered the presence of MPPS trajectories in quasilinear systems of ordinary differential equations. In the interest of brevity, MPPS trajectories are the ones which can be decomposed as the sum of periodic and Poisson stable functions. Motivated by the importance of oscillations in real world processes \cite{Cantero16, Doelling21, Gulev15, Samuelson71, Vance95} and various application fields of dynamic equations on time scales \cite{Yli14, Yli21, Li21, Liao15, Seiffertt19, Thomas09, Tisdell08}, in this study, we investigate the existence, uniqueness, and asymptotic stability of MPPS solutions in such equations. To the best of our knowledge this is the first time in the literature that Poisson stable as well as MPPS solutions are introduced and investigated for dynamic equations on time scales.

In the literature, the concept of dynamic equations on time scales has started with Hilger \cite{Hilger88}. This concept, in general, unifies the studies of differential and difference equations. The basic definitions concerning dynamic equations on time scales are as follows \cite{Agarwal02,Bohner01,Lakshmikantham96}. A time scale is a nonempty closed subset of $\mathbb R$. On a time scale $\mathbb T$, the forward jump operator is defined by $\sigma(t)=\inf \left\lbrace s\in\mathbb T ~:~ s>t\right\rbrace $, whereas $\rho(t) = \sup\left\lbrace s\in\mathbb T ~:~  s<t \right\rbrace $ is the backward jump operator. A point $t \in \mathbb T$ is called right-dense if $\sigma(t)=t$ and it is called right-scattered if $\sigma (t)>t$. Similarly, $t \in \mathbb T$ is said to be left-dense, left-scattered if $\rho(t)=t$, $\rho(t)<t$ holds, respectively. We say that a function $u:\mathbb T \to \mathbb R^m $ is rd-continuous if it is continuous at each right-dense point and its left-sided limit exists in each left-dense point. If $t$ is a right-scattered point of $\mathbb T$, then the delta derivative $u^{\Delta}$ of a continuous function $u$ is defined to be 
\begin{eqnarray}\label{deltaderivat}
u^{\Delta}(t) = \frac{u(\sigma(t)) - u(t)}{\sigma(t)-t}. 
\end{eqnarray}
Additionally, 
\begin{eqnarray} \label{deltaderivat2}
u^{\Delta}(t) = \lim_{r \to t, r \in\mathbb T} \frac{u(t) - u(r)}{t-r}
\end{eqnarray}
at a right-dense point $t$, provided that the limit exists.

It was shown by Akhmet and Turan \cite{Akhmet06} that dynamic equations on time scales which are union of disjoint closed intervals with positive length can be transformed to systems of impulsive differential equations. In the present study we make use of the technique introduced in \cite{Akhmet06} to investigate MPPS solutions of dynamic equations on a periodic time scale. More precisely, we take into account the time scale 
\begin{eqnarray} \label{timescale}
\displaystyle \mathbb{T}_0=\bigcup_{k=-\infty}^{\infty} \left[ \theta_{2k-1}, \theta_{2k} \right],
\end{eqnarray}
where for each integer $k$ the terms of the sequence $\left\lbrace \theta_k\right\rbrace_{k\in\mathbb Z}$ are defined by the equations 
\begin{eqnarray} \label{timescale11}
\theta_{2k-1}=\theta+ \delta + (k-1)\omega, \ \  \theta_{2k}=\theta + k\omega, 
\end{eqnarray}
in which $\theta$ is a fixed real number and $\omega$, $\delta$ are positive numbers such that $\omega > \delta$. 
The time scale $\mathbb T_0$ is periodic since $t \pm \omega\in \mathbb T_0$ whenever $t\in\mathbb T_0$, and 
$$\theta_{2k+1} -\theta_{2k}=\delta, \ \ k\in \mathbb Z.$$ For details of periodic time scales the reader is referred to \cite{Kaufmann06}, and some applications of dynamic equations on such time scales can be found in \cite{Bohner07,Du10,Fen17}.
It is worth noting that for each $k \in \mathbb Z$, the points $\theta_{2k}$ are right-scattered and left-dense, the points $\theta_{2k-1}$ are left-scattered and right-dense, and $\sigma(\theta_{2k})=\theta_{2k+1}$, $\rho(\theta_{2k+1}) = \theta_{2k}$.

Our main object of investigation is the equation
\begin{eqnarray} \label{maineqn}
y^{\Delta}(t)	= Ay(t) + f(t) + g(t),
\end{eqnarray}
where  $t \in \mathbb T_0$, $A \in \mathbb R^{m\times m}$ is a constant matrix,  $f:\mathbb T_0 \to \mathbb R^m$ is an rd-continuous function such that 
\begin{eqnarray} \label{periodicityf}
f(t + \omega) = f(t)
\end{eqnarray}
for each $t \in \mathbb T_0$, the function $g:\mathbb T_0 \to \mathbb R^m$ is defined by 
\begin{eqnarray} \label{funcg}
g(t)=\gamma_k
\end{eqnarray}
for $t \in [\theta_{2k-1}, \theta_{2k}]$, $k \in \mathbb Z$, and $\left\lbrace \gamma_k \right\rbrace_{k\in\mathbb Z} $ is a bounded sequence in $\mathbb R^m$.
In this paper we rigorously prove that if the sequence $\left\lbrace \gamma_k \right\rbrace_{k\in\mathbb Z} $ is positively Poisson stable, then system (\ref{maineqn}) possesses a unique asymptotically stable MPPS solution. 

It was shown in \cite{Akhmet21} that an MPPS function is Poisson stable if the corresponding Poisson number is zero. However, in this study, the structure of the function $g(t)$ in (\ref{maineqn}), which is defined by means of the sequence $\{\gamma_k\}_{k\in\mathbb Z}$, allows us to make discussion for the Poisson stability of the MPPS solution without taking into account the Poisson number.

The rest of the paper is organized as follows. In Section \ref{sec2}, we utilize the technique introduced in \cite{Akhmet06} to reduce (\ref{maineqn}) to an impulsive system. We investigate the presence of bounded solutions of the reduced impulsive system and hereby the ones for (\ref{maineqn}) in Section \ref{sec3}. The new definitions of positively Poisson stable and MPPS functions defined on time scales are provided in Section \ref{sec4}. Moreover, in that section we rigorously prove the existence and uniqueness of an asymptotically stable MPPS solution of system (\ref{maineqn}). Section \ref{sec5}, on the other hand, is devoted to an example, which confirms the theoretical results. Finally, some concluding remarks are provided in Section \ref{sec6}.

\section{Reduction to impulsive systems} \label{sec2}

In this section we make use of the $\psi$-substitution method introduced by Akhmet and Turan \cite{Akhmet06,Akhmet09} to reduce dynamic equation (\ref{maineqn}) to a system of impulsive differential equations. 

We assume without loss of generality that $\theta_{-1} <0 \leq \theta_0$. On the set $$\mathbb T'_0=\mathbb T_0 \setminus \left\lbrace \theta_{2k-1}:~ k\in\mathbb Z \right\rbrace,$$ let us consider the $\psi$-substitution defined through the equation
\begin{eqnarray} \label{func_psi}
\psi(t)= t - k \delta, \ \ \theta_{2k-1}<t\leq\theta_{2k}
\end{eqnarray}
for each integer $k$  \cite{Akhmet06}.
The function $\psi:\mathbb T'_0 \to \mathbb R$ is one-to-one and onto, $\displaystyle \lim_{t\to\infty, ~t\in \mathbb T_0'}\psi(t) = \infty$, $\psi(0)=0$, and 
\begin{eqnarray} \label{psiinv}
\psi^{-1}(s)= s + k \delta, \ \ s_{k-1} < s \leq s_{k}
\end{eqnarray}
for each integer $k$, where $s_k=\psi(\theta_{2k})$. Equation (\ref{func_psi}) yields
\begin{eqnarray} \label{defnsk}
s_k=\theta+k(\omega-\delta),
\end{eqnarray}
and accordingly, we have $s_{k+1}=s_k+\omega-\delta$ for $k\in\mathbb Z$.

The function $\psi^{-1}:\mathbb R \to \mathbb T'_0$ defined by (\ref{psiinv}) is piecewise continuous with discontinuities of the first kind at the points $s_k$, $k\in \mathbb Z$, such that $\psi^{-1}(s_k+)-\psi^{-1}(s_k)=\delta$ with $\displaystyle \psi^{-1}(s_k+)=\lim_{s\to s_k^+} \psi^{-1}(s)$. Moreover, $d\psi(t)/dt=1$ for $t\in \mathbb T'_0$ and $d\psi^{-1}(s)/ds=1$ for $s \in \mathbb R \setminus \left\lbrace s_k :~ k\in\mathbb Z \right\rbrace$ \cite{Akhmet06}. 
Since $\psi(\omega)=\omega-\delta$, one can attain by means of Corollary 12 \cite{Akhmet06} that the equality 
\begin{eqnarray} \label{psieqn1}
\psi(t+\omega)=\psi(t)+\omega-\delta
\end{eqnarray}
is fulfilled for each $t\in \mathbb T'_0$.

In what follows $C_{rd}(\mathbb T)$ stands for the set of all functions $\varphi(t):\mathbb T \to \mathbb R^m$ which are rd-continuous on a time scale $\mathbb T$.

As a consequence of Lemma 13 and Lemma 14 mentioned in paper \cite{Akhmet06}, we have the following assertion.
\begin{lemma} \label{periodiclemma} 
A function $\varphi(t) \in C_{rd}(\mathbb T_0)$ is periodic with period $\omega$ if and only if the function $\varphi(\psi^{-1}(s))$ is periodic on $\mathbb R$ with period $\omega-\delta$.  
\end{lemma}

Utilizing the descriptions of the delta derivative at right-scattered and right-dense points given respectively by (\ref{deltaderivat}) and (\ref{deltaderivat2}), one can express system (\ref{maineqn}) in the form
\begin{eqnarray} \label{system1}
	&& y'(t)=Ay(t) +f(t)+g(t),  \ \ t\in \mathbb T'_0, \nonumber \\
	&& y(\theta_{2k+1}) = y(\theta_{2k}) + \delta \left( Ay(\theta_{2k}) + f(\theta_{2k}) + \gamma_k\right) . 
\end{eqnarray}
The substitution $s=\psi(t)$, where $\psi(t)$ is defined by (\ref{func_psi}), transforms (\ref{system1}) to the impulsive system
\begin{eqnarray} \label{systemimpulse} 
&& x'(s) = Ax(s) + f\left(\psi^{-1} (s)\right) +g \left(\psi^{-1} (s)\right), \ \ s\neq s_k, \nonumber \\
&& \Delta x \big{|}_{s=s_k} =\delta \left(  A x(s_k) + f\left(\psi^{-1} (s_k)\right) + \gamma_k\right) ,
\end{eqnarray}
where $x(s) = y\left( \psi^{-1}(s)\right) $, $\Delta x \big{|}_{s=s_k}=x(s_k+)-x(s_k)$, $x(s_k+) = \displaystyle \lim_{s\to s_k^+}x(s)$, and the sequence $\{s_k\}_{k\in\mathbb Z}$ of impulse moments is defined by (\ref{defnsk}).

It is worth noting that if a function $\widetilde x(s):\mathbb R \to \mathbb R^m$ is a solution of the impulsive system (\ref{systemimpulse}), then the function $\widetilde y(t):\mathbb T_0 \to \mathbb R^m$ defined by $\widetilde y(t) = \widetilde x(\psi(t))$ for $t \in \mathbb T'_0$ with $\widetilde y(\theta_{2k+1}) =  \widetilde x(s_k+),$ $k\in\mathbb Z,$ is a solution of (\ref{maineqn}), and vice versa.

The existence and uniqueness of bounded solutions for systems (\ref{maineqn}) and (\ref{systemimpulse}) are investigated in the next section.

\section{Bounded solutions} \label{sec3}

In the remaining parts of the paper we will denote by $i(J)$ the number of the terms of the sequence $\left\{s_k\right\}_{k\in\mathbb Z}$ which take place in an interval $J$. 
One can confirm using (\ref{defnsk}) that
\begin{eqnarray} \label{equality11}
i([r+\omega-\delta,s+\omega-\delta))=i([r,s))
\end{eqnarray}
for every $s, r \in \mathbb R$ with $s > r$.

Let us denote by $U(s,r)$ the matriciant \cite{Akh1,Samolienko95} of the linear homogeneous impulsive system
\begin{eqnarray*} \label{systemimpulse2} 
&& x'(s) = Ax(s), \ \ s\neq s_k, \nonumber \\
&& \Delta x \big{|}_{s=s_k} = \delta A x(s_k)	
\end{eqnarray*}
such that $U(s,s)=I$. The equation 
\begin{eqnarray} \label{matriciant1}
U(s,r)=e^{A(s-r)} (I+\delta A)^{i([r,s))}
\end{eqnarray}
is fulfilled for $s>r$.

The following assumptions are required.
\begin{itemize}
	\item [(A1)] $\det (I+\delta A) \neq 0$, where $I$ is the $m \times m$ identity matrix;
	\item [(A2)] All eigenvalues of the matrix $e^{(\omega-\delta)A}(I+\delta A)$ lie inside the unit circle.
\end{itemize}

In the sequel we use the Euclidean norm for vectors and the spectral norm for square matrices.
Under the assumptions (A1) and (A2) there exist real numbers $N \geq 1$ and $\lambda >0$ such that 
\begin{eqnarray} \label{Cacuhyineq}
\left\| U(s,r)  \right\| \leq N e^{-\lambda (s-r)}
\end{eqnarray}
for $s \geq r$ \cite{Samolienko95}.

It is demonstrated in Theorem 87 \cite{Samolienko95} that the impulsive system (\ref{systemimpulse}) possesses a unique solution $\phi(s)$ which is bounded on the real axis and satisfies the equation
\begin{eqnarray} \label{bddsolnimpulse}
\phi(s) & =& \displaystyle \int_{-\infty}^{s} U(s,r) \left( f\left(\psi^{-1}(r)\right) + g\left(\psi^{-1}(r)\right) \right) dr \nonumber \\ 
&+& \delta \displaystyle \sum_{-\infty < s_k <s} U(s,s_{k}+) \left( f\left(\psi^{-1}(s_k)\right) + \gamma_k \right),
\end{eqnarray}
provided that (A1) and (A2) hold.
It can be verified that 
\begin{eqnarray} \label{inequality1}
\Big\| \displaystyle \int_{-\infty}^{s} U(s,r) \left( f\left(\psi^{-1}(r)\right) + g\left(\psi^{-1}(r)\right) \right) dr \Big\| \leq \frac{N(M_f + M_{\gamma})}{\lambda}
\end{eqnarray}
and
\begin{eqnarray} \label{inequality2}
\Big\|  \displaystyle \sum_{-\infty < s_k <s} U(s,s_{k}+) \left( f\left(\psi^{-1}(s_k)\right) + \gamma_k \right) \Big\| \leq \frac{N  (M_f+M_{\gamma})}{1-e^{-\lambda (\omega-\delta)}},
\end{eqnarray}
where 
\begin{eqnarray} \label{Mfdefn}
M_f=\displaystyle \sup_{t \in \mathbb T_0}  \left\|f(t)\right\|
\end{eqnarray}
and
\begin{eqnarray} \label{Mgammadefn}
M_{\gamma}=\displaystyle \sup_{k \in \mathbb Z}  \left\|\gamma_k\right\|.
\end{eqnarray}
The inequalities (\ref{inequality1}) and (\ref{inequality2}) imply that
\begin{eqnarray*}
\displaystyle \sup_{s \in \mathbb R} \left\|\phi(s)\right\| \leq N(M_f + M_{\gamma}) \left(\frac{1}{\lambda}+\frac{\delta}{1-e^{-\lambda (\omega-\delta)}}\right).
\end{eqnarray*}
Therefore, the function $\vartheta(t) : \mathbb T_{0} \to \mathbb R^m$ defined by 
\begin{eqnarray} \label{bddslntimescale1}
\vartheta(t) = \phi(\psi(t)), \ \ t \in \mathbb T'_{0},
\end{eqnarray}
and satisfying 
\begin{eqnarray} \label{bddslntimescale2}
\vartheta(\theta_{2k+1}) =  \phi(s_k+), \ \ k\in\mathbb Z,
\end{eqnarray}
is the unique solution of system (\ref{maineqn}) which is bounded on $\mathbb T_0$ such that
\begin{eqnarray*}
\displaystyle \sup_{t \in \mathbb T_0} \left\|\vartheta(t)\right\| \leq N(M_f + M_{\gamma}) \left(\frac{1}{\lambda}+\frac{\delta}{1-e^{-\lambda (\omega-\delta)}}\right).
\end{eqnarray*}

The asymptotic stability of the bounded solution $\vartheta(t)$ is discussed in the following assertion.

\begin{lemma} \label{asymptotic}
If the assumptions (A1) and (A2) are valid, then the bounded solution $\vartheta(t)$ of system (\ref{maineqn}) is asymptotically stable.
\end{lemma}

\noindent \textbf{Proof.} Let us consider a solution $\widetilde{\vartheta}(t)$ of system (\ref{maineqn}) satisfying $\widetilde{\vartheta}(t_0)=\vartheta_0$ for some $t_0\in\mathbb T'_0$ and $\vartheta_0 \in \mathbb R^m$.
We denote $\widetilde{\phi}(s) = \widetilde{\vartheta} \left( \psi^{-1}(s)\right)$ and define $\phi_0=\phi \left(\psi(t_0)\right)$. 
For $s > \psi(t_0)$, using the equations
\begin{eqnarray*}
\phi(s) & = & U\left(s,\psi(t_0)\right) \phi_0 + \displaystyle \int_{\psi(t_0)}^s U(s,r) \left( f\left(\psi^{-1}(r)\right) + g\left(\psi^{-1}(r)\right) \right) dr \\
& + & \delta \displaystyle \sum_{\psi(t_0) < s_k <s} U(s,s_{k}+) \left( f\left(\psi^{-1}(s_k)\right) + \gamma_k \right)
\end{eqnarray*}
and
\begin{eqnarray*}
\widetilde{\phi}(s) & = & U\left(s,\psi(t_0)\right) \vartheta_0 + \displaystyle \int_{\psi(t_0)}^s U(s,r) \left( f\left(\psi^{-1}(r)\right) + g\left(\psi^{-1}(r)\right) \right) dr \\
& + & \delta \displaystyle \sum_{\psi(t_0) < s_k <s} U(s,s_{k}+) \left( f\left(\psi^{-1}(s_k)\right) + \gamma_k \right),
\end{eqnarray*}
we obtain
\begin{eqnarray*}
\big\|\phi(s)-\widetilde{\phi}(s)\big\| \leq N \left\| \phi_0 - \vartheta_0\right\| e^{-\lambda(s-\psi(t_0))}.
\end{eqnarray*}
Thus,
\begin{eqnarray*}
\big\|\vartheta(t)-\widetilde{\vartheta}(t)\big\| \leq N \left\| \phi_0 - \vartheta_0\right\| e^{-\lambda(\psi(t)-\psi(t_0))}, \ \ t > t_0.
\end{eqnarray*}
The last inequality implies that the bounded solution $\vartheta(t)$ is asymptotically stable. $\square$

The main result of the present paper is provided in the next section.

\section{Modulo periodic Poisson stable solutions} \label{sec4}

The following definition is concerned with positively Poisson stable sequences \cite{Sell71}.

\begin{definition} (\cite{Sell71}) \label{defn1}
A bounded sequence $\left\lbrace \gamma_k \right\rbrace_{k\in\mathbb Z} $ in $\mathbb R^m$ is called positively Poisson stable if there exists a sequence $\left\lbrace \zeta_n\right\rbrace_{n \in \mathbb N}$ of positive integers which diverges to infinity such that $ \left\| \gamma_{k + \zeta_n} -\gamma_k\right\| \to 0$ as $n \to \infty$ for each $k$ in bounded intervals of integers.
\end{definition}

The definitions of positively Poisson stable and MPPS functions on time scales are as follows.

\begin{definition} \label{defn2}
Let $\mathbb T$ be a time scale such that $\sup \mathbb T =\infty$. A bounded function $\varphi(t) \in C_{rd}(\mathbb T)$ is called positively Poisson stable if there exists a sequence $\left\lbrace \eta_n \right\rbrace_{n\in\mathbb N}$ which diverges to infinity such that $\left\|\varphi(t + \eta_n) -\varphi(t)\right\| \to 0$ as $n \to \infty$ uniformly on compact subsets of $\mathbb T$. 
\end{definition}

\begin{definition} \label{defn3}
Let $\mathbb T$ be a time scale such that there exists a positive number $\omega$ with $t\pm\omega \in\mathbb T$ whenever $t\in\mathbb T$. A function $\varphi(t) \in C_{rd}(\mathbb T)$ is called a modulo periodic Poisson stable function if $\varphi(t)=\varphi_1(t) + \varphi_2(t)$ for every $t\in\mathbb T$ in which the function $\varphi_1 \in C_{rd}(\mathbb T)$ is periodic and $\varphi_2 \in C_{rd}(\mathbb T)$ is positively Poisson stable.
\end{definition}

The main result of the present study is mentioned in the following theorem.	
	
\begin{theorem} \label{mainthm}
Suppose that the assumptions (A1) and (A2) are fulfilled. If the sequence $\{\gamma_k\}_{k\in\mathbb Z}$ is positively Poisson stable, then system (\ref{maineqn}) possesses a unique asymptotically stable MPPS solution.
\end{theorem}	
	
\noindent{\textbf{Proof.}} The bounded solution $\phi(s)$ of system (\ref{systemimpulse}), which is defined by (\ref{bddsolnimpulse}), can be expressed in the form $$\phi(s) = \phi_1(s) + \phi_2(s), \ \ s\in\mathbb R,$$ where
\begin{eqnarray} \label{eqq2}
\phi_1(s) = \displaystyle \int_{-\infty}^s U(s,r) f\left(\psi^{-1}(r)\right) dr + \delta \sum_{-\infty < s_k <s} U(s,s_k+) f\left(\psi^{-1}(s_k)\right)
\end{eqnarray}
and
\begin{eqnarray} \label{eqq3}
\phi_2(s) = \displaystyle \int_{-\infty}^s U(s,r) g\left(\psi^{-1}(r)\right) dr + \delta \sum_{-\infty < s_k <s} U(s,s_k+) \gamma_k.
\end{eqnarray}
The bounded solution $\vartheta(t)$ of system (\ref{maineqn}), given by (\ref{bddslntimescale1}) and (\ref{bddslntimescale2}), satisfies the equation
\begin{eqnarray} \label{bddslntimescale5}
\vartheta(t) = \vartheta_1(t)+ \vartheta_2(t),
\end{eqnarray}
in which the functions $\vartheta_1(t) : \mathbb T_0 \to \mathbb R^m$ and $\vartheta_2(t) : \mathbb T_0 \to \mathbb R^m$ are respectively defined by
\begin{eqnarray*}  
\vartheta_1(t) = \phi_1(\psi(t))
\end{eqnarray*}
and
\begin{eqnarray*}  
\vartheta_2(t) = \phi_2(\psi(t))
\end{eqnarray*}
such that  the equations
$ 
\vartheta_1(\theta_{2k+1}) =  \phi_1(s_k+)
$
and 
$  
\vartheta_2(\theta_{2k+1}) =  \phi_2(s_k+)
$
are fulfilled for each $k\in\mathbb Z$.

The function $\vartheta(t)$ is an asymptotically stable solution of (\ref{maineqn}) by Lemma \ref{asymptotic}. In the proof, we will show that $\vartheta(t)$ is an MPPS function by respectively demonstrating the periodicity and Poisson stability of $\vartheta_1(t)$ and $\vartheta_2(t)$ in accordance with Definition \ref{defn3}. 

Firstly, let us discuss the periodicity of $\vartheta_1(t)$. We attain by means of (\ref{eqq2}) that
\begin{eqnarray} \label{prooff1}
\phi_1(s+\omega-\delta) & = & \displaystyle \int_{-\infty}^{s+\omega-\delta} U(s+\omega-\delta,r) f\left(\psi^{-1}(r)\right) dr \nonumber \\ 
&+& \delta \sum_{-\infty < s_k <s+\omega-\delta} U(s+\omega-\delta,s_k+) f\left(\psi^{-1}(s_k)\right) \nonumber \\
&=& \displaystyle \int_{-\infty}^{s} U(s+\omega-\delta,r+\omega-\delta) f\left(\psi^{-1}(r+\omega-\delta)\right) dr \nonumber \\
&+& \delta \sum_{-\infty < s_k <s} U(s+\omega-\delta,s_{k+1}+) f\left(\psi^{-1}(s_{k+1})\right). 
\end{eqnarray}
Equations (\ref{equality11}) and (\ref{matriciant1}) yield $U(s+\omega-\delta, r+\omega-\delta)=U(s,r)$ for $s>r$. In accordance with (\ref{periodicityf}), Lemma \ref{periodiclemma} implies that $f\left(\psi^{-1}(r+\omega-\delta)\right)=f\left(\psi^{-1}(r)\right)$ for $r\in\mathbb R$ and $f\left(\psi^{-1}(s_{k+1})\right)=f\left(\psi^{-1}(s_k)\right)$ for $k\in\mathbb Z$. Therefore, $\phi_1(s)$ is $(\omega-\delta)$-periodic on $\mathbb R$ by (\ref{prooff1}). Utilizing Lemma \ref{periodiclemma} one more time we obtain that the function $\vartheta_1(t)$ is $\omega$-periodic on $\mathbb T_0$.

Next, we will prove that $\vartheta_2(s)$ is positively Poisson stable. For that purpose, let us consider a fixed compact subset $\mathcal C$ of the time scale $\mathbb T_0$. There exist integers $\alpha$ and $\beta$ with $\beta > \alpha$ such that $\mathcal C \subseteq \left[\theta_{2\alpha}, \theta_{2\beta}\right] \cap \mathbb T_0$. Accordingly we have $\psi(\mathcal C) \subseteq \left[s_{\alpha}, s_{\beta}\right]$.

Take an arbitrary positive number $\varepsilon$ and a positive number $\tau_0$ with
\begin{eqnarray} \label{ineqrho}
\tau_0 \leq \displaystyle \frac{1}{2N(1+2M_{\gamma})} \left(\displaystyle \frac{1}{\lambda} + \frac{\delta}{1-e^{-\lambda(\omega-\delta)}}\right)^{-1}.
\end{eqnarray}  
Moreover, suppose that $j$ is a sufficiently large positive integer satisfying 
\begin{eqnarray} \label{ineqj}
j \geq \displaystyle \frac{1}{\lambda (\omega-\delta)} \ln \left(\frac{1}{\tau_0\varepsilon} \right).
\end{eqnarray}  
Because $\{\gamma_k\}_{k\in\mathbb Z}$ is positively Poisson stable, there is a sequence $\left\{\zeta_n\right\}_{n \in \mathbb N}$ of positive integers which diverges to infinity such that $\left\|\gamma_{k+\zeta_n} - \gamma_k\right\| \to 0$ as $n \to \infty$ for each $k$ in bounded intervals of integers. Thus, there exists a natural number $n_0$ such that for $n \geq n_0$ the inequality 
\begin{eqnarray} \label{ineqgamma}
\left\|\gamma_{k+\zeta_n} - \gamma_k\right\| < \tau_0 \varepsilon
\end{eqnarray} 
holds for each $k=\alpha-j+1, \alpha-j+2, \ldots, \beta$. Accordingly, if $n \geq n_0$, then the inequality
\begin{eqnarray} \label{ineqq1}
\left\|g\left(\psi^{-1}(s + \mu_n)\right)-g\left(\psi^{-1}(s)\right)\right\| < \tau_0 \varepsilon
\end{eqnarray}
holds for $s_{\alpha-j}<s\leq s_{\beta}$, where $\mu_n=(\omega-\delta) \zeta_n$, $n \in \mathbb N$.

Let us fix a natural number $n$ such that $n \geq n_0$. Making benefit of (\ref{eqq3}) one can obtain 
\begin{eqnarray*}
\phi_2(s+\mu_n) - \phi_2(s) &=& \displaystyle \int_{-\infty}^{s} U(s,r) \left( g\left(\psi^{-1}(r+\mu_n)\right) -g\left(\psi^{-1}(r)\right) \right) dr \\
&+& \delta \displaystyle \sum_{-\infty < s_k < s} U(s,s_k+) \left(\gamma_{k+\zeta_n} - \gamma_k\right).
\end{eqnarray*}
Therefore, for $s_{\alpha-j} \leq s \leq s_{\beta}$, we have
\begin{eqnarray*}
\left\|\phi_2(s+\mu_n) - \phi_2(s)\right\| &\leq & \displaystyle \int_{-\infty}^{s_{\alpha-j}} Ne^{-\lambda (s-r)}  \left\| g\left(\psi^{-1}(r+\mu_n)\right) -g\left(\psi^{-1}(r)\right)  \right\| dr \\
&+& \displaystyle \int^{s}_{s_{\alpha-j}} Ne^{-\lambda (s-r)}   \left\| g\left(\psi^{-1}(r+\mu_n)\right) -g\left(\psi^{-1}(r)\right)  \right\| dr \\
&+& \delta \displaystyle \sum_{-\infty < s_k \leq s_{\alpha-j}} Ne^{-\lambda (s-s_k)}  \left\| \gamma_{k+\zeta_n} - \gamma_k \right\| \\ 
&+& \delta \displaystyle \sum_{s_{\alpha-j} < s_k < s}  Ne^{-\lambda (s-s_k)}  \left\| \gamma_{k+\zeta_n} - \gamma_k \right\|.
\end{eqnarray*}
In compliance with (\ref{ineqgamma}) and (\ref{ineqq1}), it can be verified that 
\begin{eqnarray*}
\left\|\phi_2(s+\mu_n) - \phi_2(s)\right\| &\leq & \displaystyle 2NM_{\gamma} \int_{-\infty}^{s_{\alpha-j}} e^{-\lambda (s-r)} dr + N \tau_0 \varepsilon \displaystyle \int^{s}_{s_{\alpha-j}}  e^{-\lambda (s-r)} dr \\
&+& 2 \delta N M_{\gamma} \displaystyle \sum_{-\infty < s_k \leq s_{\alpha-j}} e^{-\lambda (s-s_k)} + \delta N \tau_0 \varepsilon \displaystyle \sum_{s_{\alpha-j} < s_k < s}  e^{-\lambda (s-s_k)} \\
&<& 2NM_{\gamma} \displaystyle \left(\frac{1}{\lambda} + \frac{\delta}{1-e^{-\lambda(\omega-\delta)}} \right) e^{-\lambda (s-s_{\alpha-j})} \\
&+& N \tau_0 \varepsilon \displaystyle \left(\frac{1}{\lambda} + \frac{\delta}{1-e^{-\lambda(\omega-\delta)}} \right)  \left(1-e^{-\lambda (s-s_{\alpha-j})}\right).
\end{eqnarray*}
For $s \geq s_{\alpha}$, the inequality $e^{-\lambda (s-s_{\alpha-j})} \leq \tau_0 \varepsilon$ is fulfilled since (\ref{ineqj}) is valid.
Hence, if $s_{\alpha} \leq s \leq s_{\beta}$, then
\begin{eqnarray*}
\left\|\phi_2(s+\mu_n) - \phi_2(s)\right\| < (1+2M_{\gamma}) \displaystyle \left(\frac{1}{\lambda} + \frac{\delta}{1-e^{-\lambda(\omega-\delta)}} \right)  N \tau_0 \varepsilon. 
\end{eqnarray*}
One can confirm using (\ref{ineqrho}) that
\begin{eqnarray} \label{ineqq6}
\left\|\phi_2(s+\mu_n) - \phi_2(s)\right\| < \frac{\varepsilon}{2}, \ \ s \in \psi(\mathcal{C}).
\end{eqnarray}

Now, let us denote 
\begin{eqnarray} \label{seqetan}
\eta_n=\omega \zeta_n
\end{eqnarray}
for each $n \in \mathbb N$. The sequence $\left\{\eta_n\right\}_{n\in\mathbb N}$ diverges to infinity since the same is true for $\left\{\zeta_n\right\}_{n\in\mathbb N}$. Equation (\ref{psieqn1}) yields $\psi(t) + \mu_n = \psi(t+\eta_n)$, $n\in\mathbb N$. Hence, according to (\ref{ineqq6}) we have 
\begin{eqnarray*} \label{ineqq7}
\left\|\vartheta_2(t+\eta_n) - \vartheta_2(t)\right\| < \frac{\varepsilon}{2}, \ \ t \in \mathcal{C} \cap \mathbb T'_0
\end{eqnarray*}
and
\begin{eqnarray*} \label{ineqq8}
\left\|\vartheta_2(\theta_{2k+1}+\eta_n) - \vartheta_2(\theta_{2k+1})\right\| \leq \frac{\varepsilon}{2}, \ \ k \in \mathbb Z.
\end{eqnarray*}
Therefore,
\begin{eqnarray} \label{ineqq9}
\sup_{t\in \mathcal{C}}\left\|\vartheta_2(t+\eta_n) - \vartheta_2(t)\right\| < \varepsilon.
\end{eqnarray}
The last inequality ensures that $\left\|\vartheta_2(t + \eta_n) -\vartheta_2(t)\right\| \to 0$ as $n \to \infty$ uniformly on compact subsets of $\mathbb T_0$. In other words, the function $\vartheta_2(t)$ is positively Poisson stable. Thus, the bounded solution $\vartheta(t)$ of (\ref{maineqn}) is an MPPS function.
$\square$

In conformity with Theorem \ref{mainthm} we have the following remark.

\begin{remark} \label{remark1}
Suppose that the conditions of Theorem \ref{mainthm} are valid. Using the equation $$\left\|\vartheta_1(t+\eta_n) - \vartheta_1(t)\right\|=0, \ \ t \in \mathbb T_0,$$ together with (\ref{ineqq9}), one can obtain for an arbitrary compact subset $\mathcal C$ of $\mathbb T_0$ and an arbitrary positive number $\varepsilon$ that
\begin{eqnarray*}  
\sup_{t\in \mathcal{C}} \left\|\vartheta(t+\eta_n) - \vartheta(t)\right\|  = \sup_{t\in \mathcal{C}} \left\|\vartheta_2(t+\eta_n) - \vartheta_2(t)\right\| < \varepsilon, \ \ n \geq n_0
\end{eqnarray*}
for some natural number $n_0$, where $\left\{\eta_n\right\}_{n \in \mathbb N}$ is the sequence defined by (\ref{seqetan}) and $\vartheta(t)$ is the bounded solution of (\ref{maineqn}), which satisfies (\ref{bddslntimescale5}). For that reason $\vartheta(t)$ is positively Poisson stable. In other words, system (\ref{maineqn}) admits a positively Poisson stable solution, which is asymptotically stable. 
\end{remark}

In the next section, an example possessing an MPPS solution is provided.

\section{An example} \label{sec5}

According to the result of Theorem 4.1 \cite{Akhmet17}, the logistic map
\begin{eqnarray} \label{logistic}
z_{k+1} = 3.9 z_k(1-z_k),
\end{eqnarray}
where $k \in \mathbb Z$, admits an orbit $\{z^*_k\}_{k\in\mathbb Z}$ inside the unit interval $[0,1]$ which is positively Poisson stable in the sense of Definition \ref{defn1}. 

Let us take into account the time scale $\displaystyle \mathbb{T}_0=\bigcup_{k=-\infty}^{\infty} \left[ \theta_{2k-1}, \theta_{2k} \right]$, where $\theta_{2k-1}=8k-4$ and $\theta_{2k}= 8k+1$ for $k \in \mathbb Z$. The equations (\ref{timescale11}) are satisfied for the time scale $\mathbb T_0$ with $\omega=8$, $\delta=3$, and $\theta=1$.

We consider the system
\begin{eqnarray} \label{example1}
&& y_1^{\Delta}(t)	= -\frac{2}{5} y_1(t)+\frac{1}{5}y_2(t) + \cos \left(\frac{\pi t}{4}\right) + g_1(t), \nonumber \\
&& y_2^{\Delta}(t)	= -\frac{1}{5}y_1(t)-\frac{2}{5}y_2(t) + \sin  \left(\frac{\pi t}{2}\right) + g_2(t),
\end{eqnarray}
where $t\in \mathbb T_0$, and the functions $g_1(t):\mathbb T_0 \to \mathbb R$ and $g_2(t):\mathbb T_0 \to \mathbb R$ are respectively defined by $g_1(t)=z^*_k$ and $g_2(t)=2z^*_k$ for $t \in [\theta_{2k-1}, \theta_{2k}]$, $k \in \mathbb Z$. It is worth noting that the sequence $\{\gamma_k\}_{k\in\mathbb Z}$ given by $\gamma_k=(z^*_k, 2z^*_k)^T$ is positively Poisson stable according to Theorem 3.2 \cite{Akhmet18}. System (\ref{example1}) is in the form of (\ref{maineqn}) with
$$y(t)=(y_1(t),y_2(t))^T, \ \ A=\begin{pmatrix} \displaystyle - 2/5 & \displaystyle 1/5 \\ \displaystyle -1/5 & \displaystyle -2/5 \end{pmatrix},$$ 
$$f(t)= \left(\cos \left(\frac{\pi t}{4}\right), \sin  \left(\frac{\pi t}{2}\right) \right)^T, \ \ g(t)=(g_1(t),g_2(t))^T.$$
The matrix $e^{5A}(I+3A)$, where $I$ is the $2 \times 2$ identity matrix, admits a pair of complex conjugate eigenvalues both of which are inside the unit circle, and $\det(I+3A)=2/5$. The assumptions (A1) and (A2) are satisfied for system (\ref{example1}), and therefore, it admits a unique asymptotically stable MPPS solution by Theorem \ref{mainthm}. Moreover, the MPPS solution is at the same time positively Poisson stable according to Remark \ref{remark1}.

\section{Conclusion} \label{sec6}
We take into account a periodic time scale which is the union of infinitely many disjoint compact intervals with a positive length, and investigate the existence, uniqueness as well as asymptotic stability of MPPS solutions for dynamic equation on such time scales. In our discussions we make use of the reduction technique to impulsive systems introduced in \cite{Akhmet06}. The Poisson stability in (\ref{maineqn}) is inherited from the sequence $\left\{\gamma_k\right\}_{k \in \mathbb Z}$. The descriptions of positively Poisson stable and MPPS functions on time scales are newly introduced in the present study. Moreover, it is shown that the obtained MPPS solutions are at the same time positively Poisson stable. Even though in general MPPS functions are not necessarily positively Poisson stable \cite{Akhmet21}, this is true in our case owing to the commensurability of the periods of the time scale $\mathbb T_0$ and the function $f(t)$ used in (\ref{maineqn}). In the future, our results can be developed for differential equations on variable time scales \cite{Akhmet09}.

\end{document}